\documentstyle[12pt]{article}
\newcommand{\lab}[1]{\label{#1}}

\newcommand{\rf}[1]{(\ref{#1})}

\newcommand\beq{\begin{equation}}
\newcommand\eeq{\end{equation}}
\newcommand\beqy{\begin{eqnarray}}
\newcommand\eeqy{\end{eqnarray}}
\newcommand\beqys{\begin{eqnarray*}}
\newcommand\eeqys{\end{eqnarray*}}

\newtheorem{Th}{Theorem}
\newtheorem{Cor}{Corollary}

\newcommand{\reali}{{\hbox{{\rm I}\kern-.2em\hbox{\rm R}}}}
\newcommand{\complessi}{{\ \hbox{{\rm I}\kern-.6em\hbox{\bf C}}}}
\newcommand\de{\partial}
\newcommand{\ad}{{\,\hbox{\rm ad}}}

\newcommand{\Hol}{{\,\hbox{\rm Hol}\,}}
\newcommand{\trqw}{{\rm Tr_{\frac12}\,}}
\newcommand{\tr}[1]{{\rm Tr_{#1}\,}}

\newcommand{\xnot}{{x_0}}
\newcommand{\ave}[1]{\left\langle{#1}\right\rangle}
\newcommand{\aveBF}[2]{\ave{#1}_{BF,#2}}

\newcommand{\sbf}{S_{BF}}
\newcommand{\tsbf}{\widetilde{S_{BF}}}

\newcommand{\dodo}{{\cdot\cdot}}
\newcommand{\ftilde}{{\widetilde f}}
\newcommand{\btilde}{{\widetilde b}}
\newcommand{\Etilde}{{\widetilde E}}
\newcommand{\Jtilde}{{\widetilde J}}
\newcommand{\bino}[2]{
\left({
\begin{array}{c}
{#1}\\
{#2}
\end{array}
}\right)}

\def\cmp#1,{{ Commun.\ Math.\ Phys.\ \bf #1},}
\def\jmp#1,{{ J.\ Math.\ Phys.\ \bf #1},}
\def\pl#1,{{ Phys.\ Lett.\ \bf #1},}
\def\np#1,{{ Nucl.\ Phys.\ \bf #1},}
\def\mpl#1,{{ Mod.\ Phys.\ Lett.\ \bf #1},}
\def\pr#1,{{ Phys.\ Rev.\ \bf #1},}
\def\prl#1,{{ Phys.\ Rev.\ Lett.\ \bf #1},}
\def\lmp#1,{{ Lett.\ Math.\ Phys.\ \bf #1},}

\begin{document}
\begin{titlepage}
\title{Cabled Wilson Loops in $BF$ Theories}
\author{
Alberto S. Cattaneo\\[10pt]
Lyman Laboratory of Physics\\
Harvard University\\
{\sc Cambridge, MA 02138, USA}\\
E-mail: {\tt cattaneo@math.harvard.edu}
}
\date{February 6, 1996\\
{\sf HUTMP 96 B349, q-alg/9602015}\\[10pt]
P.A.C.S. 02.40, 11.15, 04.60}

\maketitle

\begin{abstract}
A generating function for cabled Wilson loops in 
three-dimensional $BF$ theories is defined, and a careful study of
its behavior for vanishing cosmological constant is performed.
This allows an exhaustive description of the
unframed knot invariants coming from the
pure $BF$ theory based on $SU(2)$, and
in particular, it proves a conjecture relating them to the 
Alexander--Conway polynomial.
\end{abstract}
\thispagestyle{empty}
\end{titlepage}
\setcounter{page}{2}
\section{Introduction}
In \cite{CCM}, a relation was conjectured between 
the Alexander--Conway polynomial and the v.e.v.\ 
({\em vacuum expectation value}) of a newly proposed observable
for the TQFT ({\em topological quantum field theory}) known as $BF$
theory \cite{Schw} (or as ``pure" $BF$ theory, to distinguish it
from its generalization ``with a cosmological constant").
The main purpose of this paper is to clarify that conjecture and to
prove it in the case of $SU(2)$, although in a slightly different form.
To do so, however, we also perform a careful study of $BF$ theories, both
pure and with a cosmological constant, introducing new observables
that describe cabled Wilson loops.

$BF$ theories are interesting from both the mathematical and the
physical point of view: They are TQFTs that, in principle,  can be
defined in any dimension (s.\ \cite{CCFM} and references therein). 
They describe quantum gravity in three dimensions \cite{WittGrav}
and could be useful for the study of four-dimensional quantum gravity 
in the approach of \cite{ASR}. Moreover, they are related to 
Yang--Mills theories in the weak coupling limit, where
they could be used  to describe a ``topological phase'' \cite{CCGM}
(s.\ also \cite{Ans} for a related approach). 
Finally, hitherto only in the Abelian case,
they have been used to study many-body systems \cite{FGM}. 

On the other side,
the interest for the Alexander--Conway polynomial comes from the fact
that it has always played a different role than
the other knot invariants, being the only one that, hitherto, has
been described by the methods of the traditional algebraic 
topology; yet it lacked a TQFT description (and \cite{CCM} was an
attempt to fill this gap). 

Indeed,
the Chern--Simons theory describes the knot invariants related
to the $q$-deformation of a classical group $G$ in terms
of v.e.v.'s of Wilson loops (i.e.,
traces of $G$-holonomies along the knots) \cite{Witt89}; yet the 
$q=1$ case, i.e., the Alexander--Conway polynomial, does not correspond
to any acceptable value of the Chern--Simons coupling constant $k$. 

In the case of $SU(2)$, an important improvement was given by 
the Melvin--Morton conjecture \cite{MM} (proved in \cite{BG}),
which states that the inverse of the Alexander--Conway 
polynomial appears in the $q-1$ expansion of the colored Jones function.
Rozansky \cite{Roz} was then able to prove the same result in the framework
of the Chern--Simons theory; more precisely, he showed that the inverse of
the Alexander--Conway polynomial is recovered by a saddle-point computation in
the limit $k\to\infty$.

In \cite{CCFM}, an equivalence between the v.e.v.'s of the
Chern--Simons theory and of the $BF$ theory with cosmological constant 
$\kappa=1/(2k)$ was shown.  The $k\to\infty$ limit (i.e., the $q=1$
case) is then simply described by removing the cosmological term, thus
obtaining the pure $BF$ theory.

In this paper we show that the trace in the fundamental
representation of the pure-$BF$-theory
observable \rf{B.7}, first introduced in \cite{CCM},
is related to the $\kappa\to0$ limit of a newly proposed observable,
\rf{B.3},
for the $BF$ theory with a cosmological constant. 
In the framework of the Chern--Simons
theory, this observable corresponds to the ``exponential" of a 
holonomy; so its trace is the generating function for cabled Wilson
loops, and the v.e.v.\ of this trace
is the generating function for cabled knot invariants (as suggested
by the results of \cite{KSA}). 

Since the colored Jones functions are related to the cabled knot invariants,
we are able to connect the $q-1$ expansion 
of the former to the pure-$BF$-theory knot
invariant.  By the Melvin--Morton conjecture, we conclude that
the latter is 
related to the Alexander--Conway polynomial, the precise relation
being \rf{D.23}.
Our proof, however, can also go the other way; i.e., should one 
directly prove a relation between the pure $BF$ theory and
the Alexander--Conway polynomial, then we would have a further proof
of the Melvin--Morton conjecture.

Our results are not limited to the trace in the fundamental representation 
of a particular observable for pure $BF$ theory. In fact, we are able
to prove that, in the standard framing, this observable is the most
general one can consider, s.\ Thm.\ \ref{th-sf} and Cor.\ \ref{cor-mg}.
Moreover, we show how to compute the higher-representation knot invariants,
s.\ \rf{G.5.a} and \rf{G.5.b}, and also evaluate the limit as the dimension
of the representation goes to infinity, s.\ \rf{WWtilde} and
\rf{RW}.  Thus, we have a complete description of the unframed 
knot invariants coming from the pure $BF$ theory based on $SU(2)$, 
s.\ Thm.\ \ref{th-BF}. The case of links has not been considered yet. 
\vskip 10pt

The paper is organized as follows: In Sec.\ \ref{sec-CJF}, we describe
the colored Jones functions and introduce their generating function.
In Sec.\ \ref{sec-BTCJF}, we recall the property of the $BF$ theory with a
cosmological constant and define a generating function for cabled
Wilson loops.
In Sec.\ \ref{sec-PBF}, we discuss the pure $BF$ theory and prove 
Thm.\ \ref{th-sf} and Cor.\ \ref{cor-mg}.
In Sec.\ \ref{sec-TOCC}, we consider the limit for vanishing coupling
constant and prove Thm.\ \ref{th-BF}.
We conclude our work in Sec.\ \ref{sec-PEW} comparing our present results
with those obtained perturbatively in \cite{CCM}. 
For the sake of clarity, we have put all the cumbersome computations into 
some appendixes, to which we shall refer in the text.

\section{The Colored Jones Function}
\lab{sec-CJF}
The colored Jones function, $J_d(C;h)$ is defined in \cite{MM} as the invariant
of the knot $C$
obtained using the irreducible $SU(2)_q$-module ($q=e^h$) of dimension $d$.
The main properties
of its expansion as a rational power series
(here and in the following, the notation of \cite{BG} is used),
\beq
J_d(C;h) = 
\sum_{m=0}^\infty \Jtilde_m(C;d)\ h^m=
d \sum_{j,m=0}^\infty b_{jm}(C)\ (d-1)^jh^m,
\lab{expand}
\eeq
were stated in \cite{MM}:
\begin{enumerate}
\item The functions 
$\Jtilde_m(C;d)$ are odd in $d$. This means that (the analytic
continuation of) $J_d(C;h)$ itself is odd
\beq
J_{-d}(C;h)=-J_d(C;h).
\lab{Jodd}
\eeq
(This of course implies that the coefficients $b_{jm}$ in \rf{expand} are
not completely independent.)
\item
These functions are actually polynomials of
degree not exceeding $2m+1$.  
\item In the hypothesis that 
{\em the standard framing for $C$ is chosen}, the degree of $\Jtilde_m$ 
cannot exceed $2m-1$.  
\end{enumerate}
In the same hypothesis of standard framing, two conjectures were also
stated in \cite{MM}: 
\begin{enumerate}
\item The degree of $\Jtilde_m$ cannot actually exceed $m+1$, or
in other words,
$b_\dodo$ in \rf{expand} is an upper triangular matrix,
\beq
b_{jm}(C) = 0, \quad\mbox{if $j>m$}.
\lab{Con1MM}
\eeq
\item The diagonal elements of $b_\dodo$ are related to
the Alexander--Conway polynomial;
viz., if one defines the Melvin--Morton
function as
\beq
JJ(C;\hbar) = \sum_{m=0}^\infty b_{mm}(C)\ \hbar^n=
\lim_{d\to\infty,\ h\to0,\ dh=\hbar} \frac {J_d(C;h)}d
\lab{MMf}
\eeq
(the limit existing because of \rf{Con1MM}), then 
\beq
\hbar\, JJ(C;\hbar) = \frac z{\Delta(C;z)},
\quad z=e^{\hbar/2}-e^{-\hbar/2},
\lab{Con2MM}
\eeq
where $\Delta(C;z)$ is the Alexander--Conway polynomial with the 
normalization $\Delta(\bigcirc;z)=1$, where $\bigcirc$ is the 
unknot, and skein relation
\[
\Delta(C_+;z)-\Delta(C_-;z)=z\,\Delta(C_0;z).
\]
\end{enumerate}

These two conjectures have been given a
functional-integral proof by Rozansky \cite{Roz}. In particular
he proved
that the Melvin--Morton function, \rf{MMf}, is related to the 
large-$k$ limit of the Chern--Simons theory, so it can be computed in
saddle-point approximation.

Eventually, Bar-Natan and Garoufalidis \cite{BG} gave a 
mathematical proof of \rf{Con1MM} and \rf{Con2MM}
on the level of weight systems.

\subsection*{The generating function for colored Jones functions}
In the following sections, the generating function
\beq
f(C;x,h) := \sum_{d=1}^\infty \frac{x^d}{d!} J_d(C;h),
\lab{gf}
\eeq
with $x\in\complessi$, will be needed; so it is useful to anticipate here
some of its properties. These are better clarified if one introduces 
$\ftilde$ defined by 
\beq
f(C;x,h) = x\, e^x\, \ftilde(C;x,h).
\lab{C.13}
\eeq 
In fact, the expansion of $\ftilde$ as a rational power series,
\beq
\ftilde(C;x,h) = \sum_{j,m=0}^\infty \btilde_{jm}(C)\ x^jh^m,
\lab{C.18}
\eeq
shares the same properties as the expansion \rf{expand} of $J_d(C;h)$;
viz., s.\ App.\ \ref{app-f}, 
one can prove that 
\beq
b_\dodo\ \mbox{upper triangular}
\Longleftrightarrow
\btilde_{\dodo}\ \mbox{upper triangular}.
\lab{bbtilde}
\eeq
Thus, \rf{Con1MM} implies that
$\btilde_\dodo$ is actually an upper triangular matrix,
\beq
\btilde_{jm}(C) = 0, \quad\mbox{if $j>m$},
\lab{tildeb}
\eeq
and that $\ftilde$ has a well-defined limit for $|x|\to\infty$,
$h\to 0$ and $xh=\hbar$ kept fixed.  In App.\ \ref{app-f}, it is
shown that
\beq
\lim_{|x|\to\infty,\ h\to0,\ xh=\hbar}
\ftilde(C;x,h) = JJ(C;\hbar),
\lab{C.14}
\eeq
and that, in general, for any $n\ge 0$,
\beq
\lim_{|x|\to\infty,\ h\to0,\ xh=\hbar}
\left({x\frac\de{\de x}}\right)^n \ftilde(C;x,h) 
=
\left({\hbar\frac d{d\hbar}}\right)^n JJ(C;\hbar).
\lab{C.16}
\eeq

\section{The $BF$ Theory with a Cosmological Constant}
\lab{sec-BTCJF}
{}From the field-theoretical point of view,
the colored Jones function is defined as the knot invariant obtained from the
v.e.v.\ of a suitable observable, the 
Wilson loop, in the Chern--Simons theory \cite{Witt89} (actually, the
general case of a compact simple Lie group is discussed there).

For the aims of this paper, however, it is more convenient to consider 
the formulation given in terms of the $BF$ theory with a cosmological 
constant, which is equivalent \cite{C,CCFM} to the previous one.

We will consider here only the case of a single knot $C$
imbedded in $S^3$ and study its invariants associated to the group
$SU(2)$; viz., we consider the $SU(2)$-principal bundle
$P\longrightarrow S^3$, and define the action
\beq
\tsbf(\kappa) = \frac1{2\pi}\int_{S^3} \tr{} \left({B\wedge F+
\frac{\kappa^2}3 B\wedge B\wedge B }\right),
\lab{tsbf}
\eeq
where $F$ is the curvature two-form of the connection $A$, $B$
is a form in $\Omega^1(S^3,\ad P)$, and the parameter $\kappa$ is called
the {\em cosmological constant} (this name comes from the quantum-gravity
interpretation of \rf{tsbf} \cite{WittGrav}).

Then, given a knot $C$ and a base point $\xnot\in C$, 
we consider the observable
\beq
\Gamma(C,\xnot;\kappa) := \Hol_\xnot(A+\kappa B;C) =
\sum_{n=0}^\infty \kappa^n\gamma_n(C,\xnot),
\lab{Gamma}
\eeq
where $\kappa$ is the same cosmological constant as 
in \rf{tsbf}, $\Hol_\xnot(A+\kappa B;C)$ denotes the holonomy of 
the connection $A+\kappa B$ along the knot $C$ with base point $\xnot$, 
and the $\gamma$'s are
functionals of $A$ and $B$ obtained by Taylor expanding the previous 
holonomy.

Notice that
the holonomy of a connection $A$ along a curve $C$, open or closed, 
can be written as
\beq
\Hol(A;C) = P\, \exp\left({
\int_C A
}\right),
\lab{hol}
\eeq
where $P$ denotes path-ordering. Thus, by \rf{Gamma} and \rf{hol}, 
the $\gamma$'s turn out to be \cite{C,CCFM}
iterated Chen integrals of the form
\beq
\gamma_n(C,\xnot) :=\oint_{(C,\xnot)}  
\underbrace{
\hat B\cdots\hat B}_{\mbox{$n$ times}}\ 
\Hol_{\xnot}(A;C),
\lab{chen}
\eeq
where
the iterated integral $\displaystyle{
\int_a^b \omega_1\cdot\omega_2
\cdots \omega_n}$ of $n$ one-forms $\{\omega_i\}_{i=1,\cdots n}$ 
is given by the formula 
$\displaystyle{\int_{a<x_1<\cdots<x_n<b}
\omega_1(x_1)\wedge \omega_2(x_2)\wedge\cdots\wedge\omega_n(x_n)}$.
The $su(2)$-valued one-form $\hat B$ is defined as
\beq
\hat B(x) := \Hol_{\xnot}^x B(x) [\Hol_{\xnot}^x]^{-1},
\lab{bhat}
\eeq
with the holonomies computed along the portions of the
knot $C$ going from the base point $\xnot$ to the running points $x$.

The obervable $\Gamma$ in \rf{Gamma} is invariant under a gauge transformation
of the connection $A-\kappa B$, while, under the gauge transformation
\beq
A+\kappa B \rightarrow g\, (A+\kappa B)\, g^{-1} + gdg^{-1},
\lab{gauge}
\eeq
one has
\beq
\Gamma(C,\xnot;\kappa) \rightarrow g(\xnot)\, \Gamma(C,\xnot;\kappa)
g(\xnot)^{-1}.
\lab{Gammagauge}
\eeq
Thus, taking the trace of $\Gamma$
gives a gauge invariant observable that is,
moreover, independent of $\xnot$.  The knot invariant given by the v.e.v. of
this observable, w.r.t.\ the (normalized)
Gibbs weight $\exp[i\tsbf(\kappa)]$,
can be recognized by exploiting the relation with the 
Chern--Simons theory \cite{CCFM}, viz., 
\beq
J_d(C;h) = \aveBF{\tr{(d-1)/2} \Gamma(C,\xnot;\kappa)}\kappa,
\quad h=4\pi i\kappa,
\lab{B.2}
\eeq
where $\tr s$ is the trace in the irreducible representation of spin $s$.
In this framework, \rf{Jodd} is an immediate consequence of \rf{trs}.

Notice, eventually, that the observable $\Gamma$ can be decomposed into
its even ($\Gamma_0$) and odd 
($\kappa\Gamma_1$) parts, 
\beq
\Gamma(C,\xnot;\kappa) = \Gamma_0(C,\xnot;\kappa) +
\kappa \Gamma_1(C,\xnot;\kappa),
\lab{Gamma01}
\eeq
and that both of them are good observables
for the $BF$ theory with a cosmological constant \cite{C}. 

\subsection*{The $BF$-theory generating function}
The behavior \rf{Gammagauge} of the observable $\Gamma$ under a gauge
transformation of the connection $A+\kappa B$ shows that all of the traces
of powers of $\Gamma$ are gauge invariant;  so we can consider
the knot invariants given by their v.e.v.'s, or as will be clear in 
Sec.\ \ref{sec-TOCC}, their generating function, 
\beq
\begin{array}{rcl}
E(C;x,h) 
&:=& 
\sum_{n=0}^\infty\frac{x^n}{n!}\aveBF{\trqw\Gamma(C,\xnot;\kappa)^n}\kappa
=\aveBF{\trqw\exp[x\Gamma(C,\xnot;\kappa)]}\kappa,\\
h &=& 4\pi i\kappa,
\end{array}
\lab{B.3}
\eeq
with $x\in\complessi$. We consider only the fundamental representation here,
for the discussion in App.\ \ref{app-A} shows that considering any other
representation does not give further content of information.
Notice that the rightmost term in \rf{B.3} is only formal since no addition is
defined in the group.  However, this notation usefully reminds us that 
$E$ has the formal properties of an exponential.
Notice that, since $\Gamma$ is a holonomy, one has
\beq
\Gamma(C,\xnot;\kappa)^n = \Gamma(nC,\xnot;\kappa),
\lab{Gamman}
\eeq
where by $nC$ we denote the $n$th cabling of the knot $C$; thus, 
$E(C;x,h)$ is {\em the generating function for cabled knot 
invariants}.  
In the framework of the Chern--Simons theory, $\Gamma$ is replaced by
the holonomy of the connection $A$; thus $E$ can also be seen as
{\em the generating function for cabled Wilson loops}.

By using the formulae in App.\ \ref{app-A},
one can show that there is a simple relation with the generating
function for colored Jones functions defined in \rf{gf}. Actually, by
\rf{A.7}, one can express $\trqw\Gamma^n$ in terms of the traces of $\Gamma$
in the representations of spin $\frac n2$ and $\frac n2-1$.  By \rf{B.2},
it then follows that
\beq
E(C;x,h)=\sum_{n=0}^\infty \frac{x^n}{n!}[J_{n+1}(C;h)-J_{n-1}(C;h)],
\lab{B.10}
\eeq
where, according to \rf{trs}, $J_1=1$, $J_0=0$, $J_{-1}=-1$.

By \rf{gf} and \rf{B.10}, and also noting that
$E(C;0,h)=2$ and $\de_xf(C;0,h)=J_1(C;h)=1$, one obtains
\beq 
E(C;x,h) = 1 + \de_xf(C;x,h) - \int_0^x d\xi\ f(C;\xi,h).
\lab{B.16}
\eeq

\section{The Pure $BF$ Theory}
\lab{sec-PBF}
The $BF$ theory has the nice property that its limit for vanishing cosmological
constant (corresponding to the $k\to\infty$ limit in the Chern--Simons theory)
is still represented by a TQFT known as pure $BF$ theory \cite{Schw, BT},
whose action reads
\beq
\sbf = \frac1{2\pi}\int_{S^3} \tr{} \left({B\wedge F}\right).
\lab{sbf}
\eeq
It is immediately seen that $\sbf$ is invariant under gauge transformations 
\beq
\begin{array}{rcl}
A &\rightarrow& g\, A\, g^{-1} + gdg^{-1},\\
B &\rightarrow& g\, B\, g^{-1},
\end{array}
\lab{gaugeBF}
\eeq
as well as under $B$-transformations
\beq
\begin{array}{rcl}
A &\rightarrow& A,\\
B &\rightarrow& B + d_A\psi.
\end{array}
\lab{BBF}
\eeq
Here $g$ is a map from $S^3$ to the group, while $\psi$ is a form
in $\Omega^0(S^3, \ad P)$. The action \rf{sbf} is invariant under
\rf{BBF} owing to the Bianchi identity.

The good observables for the pure $BF$ theory are of course obtained
by $\Gamma_0(C,\xnot;\kappa)$ and $\Gamma_1(C,\xnot;\kappa)$, defined in
\rf{Gamma01}, in the limit $\kappa\to0$, 
and are simply given by
$\gamma_0(C,\xnot)$ and $\gamma_1(C,\xnot)$ in \rf{chen}; viz.,
\beq
\begin{array}{rcl}
\gamma_0(C,\xnot) &=& \Hol_\xnot(A;C),\\
\gamma_1(C,\xnot) &=& \oint_{x\in (C,\xnot)} \!\!
\Hol_\xnot^x\, B(x)\,\Hol_x^\xnot.
\end{array}
\lab{gamma01}
\eeq
Under \rf{gaugeBF} and \rf{BBF}, they transform as
\beq
\gamma_i(C,\xnot) \rightarrow g(\xnot)\, \gamma_i(C,\xnot)\, 
g(\xnot)^{-1},\quad i=0,1,
\lab{transfgammas}
\eeq
provided we fix $\psi(\xnot)=0$, which can be always done \cite{CCM}.

Thus, we are led to consider 
the v.e.v.'s of $\tr d\gamma_0$ and $\tr d\gamma_1$. They, however,
are rather trivial, for they correspond, respectively, 
to $J_d(C;0)=d$ and $(d/dh)J_d(C;0)=0$.  

{}From a field-theoretical point of view, the reason of this triviality
is the following:  In perturbative $BF$ theory (on $S^3$),
non vanishing v.e.v.'s contain a number $n_A$ of fields $A$ not exceeding the 
number $n_B$ of fields $B$ \cite{CCM}.  Therefore, after expanding the
holonomies in \rf{gamma01} in powers of $A$, one sees that
only the zeroth-order term $\tr{(d-1)/2}I=d$ (where $I$ is the group identity)
survives in the v.e.v.\ of $\tr{(d-1)/2}\gamma_0$; 
while, in the v.e.v.\ of $\tr{(d-1)/2}\gamma_1$, only the first-order term
does (the zeroth-order term vanishing because it is the trace of an element of 
$su(2)$). But the latter v.e.v.\ (of the form $\aveBF{\oint_C A\oint_C B}0$) 
gives the self-linking number of the knot $C$ that, by the hypothesis
of standard framing, is zero. (If one uses a more general framing,
the self-linking number is the only information one gets from these v.e.v.'s.)

This is the reason why in \cite{CCM} one looked for composite observables
that contain a higher number $n_B$ of fields $B$, and are invariant under
\rf{gaugeBF} and \rf{BBF}. 
In particular, the generating function of the v.e.v.'s of
$\tr{(d-1)/2}(\gamma_1)^n$  was taken into account,
\beq
W_d(C;\lambda) := 
\sum_{n=0}^\infty \frac{\lambda^n}{n!}\aveBF{\tr{(d-1)/2}
\gamma_1(C,\xnot)^n}0 =
\aveBF{\tr{(d-1)/2} \exp[\lambda\gamma_1(C,\xnot)]}0,
\lab{B.7}
\eeq  
with $\lambda\in\complessi$. Again, as in \rf{B.3}, the rightmost term is only
formal.

In Sec.\ \ref{sec-TOCC}, we shall need \rf{B.7} to study the limit
$\kappa\to 0$ in \rf{B.3}, to which \rf{B.7} is apparently related.
However, there is a second, perhaps more natural, generating function
to be considered, viz.,
\beq
\widetilde W_d(C;\lambda) := \aveBF{\tr{(d-1)/2} 
\exp[\lambda\beta(C,\xnot)]}0,
\lab{Wtilde}
\eeq
where
\beq
\beta(C,\xnot):= \gamma_1(C,\xnot)\cdot\gamma_0(C,\xnot)^{-1}=
\oint_{C,\xnot} \!\!\hat B \in su(2). 
\lab{beta}
\eeq
Notice that the exponential in \rf{Wtilde} is not formal, but it actually
represents the exponential map from the algebra to the group.

The generating function $\widetilde W$ is particularly useful
because it allows to understand the meaning of the knot invariants related 
to higher-dimensional representations; 
indeed, by using \rf{gn}, one can easily prove that
\beqy
\widetilde W_{2k}(C;\lambda) &=& \sum_{l=1}^k 
\widetilde W_2(C;(2l-1)\lambda),
\lab{G.5.a}\\
\widetilde W_{2k+1}(C;\lambda) &=& 1 + \sum_{l=1}^k 
\widetilde W_2(C;2l\lambda).
\lab{G.5.b}
\eeqy
Then, a formal resummation of the geometric sums appearing in \rf{G.5.a} and
\rf{G.5.b} gives
\beq
\widetilde W_d(C;\lambda) =\left\{
\begin{array}{ll}
\aveBF{\trqw\frac{\exp[(d+1)\lambda\beta(C,\xnot)]-1}
{\exp[2\lambda\beta(C,\xnot)]-1}}0 & \mbox{if $d$ is even}\\
\aveBF{\trqw\frac{\exp[(d+1)\lambda\beta(C,\xnot)]-1}
{\exp[2\lambda\beta(C,\xnot)]-1}}0 +1 & \mbox{if $d$ is odd}
\end{array}\right.
\lab{di.1}
\eeq
These formulae show that the limit
\beq
\widetilde{WW}(C;\mu) := \lim_{d\to\infty,\lambda\to0, d\lambda = \mu}
\frac{\widetilde W_d(C;\lambda)}d
\lab{WWtilde}
\eeq
is well defined, and that
\beq
\widetilde{WW}(C;\mu) = \aveBF{\trqw\frac{\exp[\mu\beta(C,\xnot)]-1}
{2\mu\beta(C,\xnot)}}0 =
\sum_{n=0}^\infty \frac1{2(n+1)!}\aveBF{\trqw[\mu\beta(C,\xnot]^n
}0.
\lab{di.4}
\eeq
Thus, by \rf{Wtilde} and \rf{di.4}, we get
\beq
\widetilde W_2(C;\lambda) = 2\frac d{d\lambda}
[\lambda\ \widetilde{WW}(C;\lambda)],
\lab{D.23bistilde}
\eeq
or equivalently,
\beq
\widetilde{WW}(C;\mu) =\frac1\mu \int_0^{\mu/2}dx
\ \widetilde W_2(C;2x).
\lab{RW}
\eeq
Notice that \rf{G.5.a} and \rf{G.5.b} are actually discretizations
of the Riemann integral in \rf{RW} with $x=l\lambda$.

A more rigorous proof of the above results can be obtained by expanding
$\widetilde W_2$ in powers of $\lambda$ and using \rf{G.5.a} and \rf{G.5.b}
to get the corresponding expansion of $\widetilde W_d$.  Then one can
use the asymptotic formula
\[
\sum_{l=1}^k l^n \stackrel{k\to\infty}\sim \frac{k^{n+1}}{n+1}
\]
to get the rightmost term in \rf{di.4}.
\vskip 10pt

So far we have considered two different generating functions, viz.,
\rf{B.7} and \rf{Wtilde}; in general, however, there are even more choices. 
Indeed, because of \rf{transfgammas},
a product of powers of
$\gamma_0$ and $\gamma_1$ (in any order) is still an observable whose trace is
invariant under \rf{gaugeBF} and \rf{BBF}.

To get rid of this arbitrariness in the definition
of the generating function, we need the following
\begin{Th}
In pure $BF$ theory, the knot invariant obtained by the v.e.v.\ of the trace
of a function of $\gamma_0(C,\xnot)$ and $\gamma_1(C,\xnot)$ is equal
to the knot invariant obtained by replacing
$\gamma_0(C,\xnot)$ with the identity of the group, {\sf
provided the standard framing is chosen}.
\lab{th-sf}
\end{Th}

In \cite{C}, a proof of the theorem 
based on an explicit study of the Feynman integrals appearing in the 
evaluations of the v.e.v.'s was given.  Here, however, we give a
simpler argument, based on the formal properties of the 
pure $BF$ theory.  

First of all we remind that, in the explicit
evaluation of the v.e.v.'s in $BF$ theory, the choice of framing 
is done by evaluating all
the holonomies in \rf{gamma01} on a companion knot $C'$ obtained by $C$
in terms of a ``small,'' non-vanishing normal displacement
\cite{CCM}. Notice that the self-linking
number of $C$ is then, by definition, the linking number of $C$ and $C'$. 

The proof of the theorem, then, essentially relies on the fact 
that, if this linking number vanishes (standard framing),
one can deform $C'$ in
$\gamma_0$ in such a way that it can be completely unlinked from $C$ and,
as such, shrunk to a point; so its holonomy, $\gamma_0$, becomes the identity.

This is possible since a knot $C'$ appearing in $\gamma_0$ (not necessarily
a framing for $C$) ``does not see
itself.''\@ In fact, a small deformation of the knot $C'$ at a
certain point $x$ amounts to introducing curvature terms in $x$ (remember that
$\gamma_0$ is a holonomy). But in pure $BF$ theory, the curvature vanishes
everywhere but on the knot $C$, along which $\gamma_1$ is evaluated.
(We refer to \cite{CCM} for the proof that, in an observable,
$A$ acts as a source for $d_AB$, while $B$ acts
as a source for the curvature $F$.)\@ Therefore, $C'$ can be freely deformed
as far as this deformation does not intersect $C$; self-intersections 
of $C'$ are allowed, however.

When such a self-intersection occurs, one can split 
$C'$ into two closed curves $C'_1$ and $C'_2$; correspondingly,
the holonomy $\gamma_0(C')$ can be written as 
$\gamma_0(C'_1)\cdot\gamma_0(C'_2)$,

To prove the theorem, one just has to repeat
this procedure until $C'$ is replaced by a collection
of circles $C'_i$ all of which surround one single strand of $C$. Then one
moves the $C'_i$'s along $C$, and cuts and splices them together 
again to form a circle $C''$.
 
Since all these deformations do not affect the linking between $C$
and $C'$, the linking number
of $C$ and $C''$ is the same as the linking number of
$C$ and $C'$. If this linking number vanishes (standard framing), then 
the circle $C''$ and the knot $C$ are unlinked, and $C''$ can be srhrunk
to a point. This concludes the proof of Thm.\ \ref{th-sf}.

As a consequence of Thm.\ \ref{th-sf}, we have the following
\begin{Cor}
$W_d(C;\lambda)$ is the most general v.e.v.\ one
can consider in pure $BF$ theory if the standard framing is chosen.
\lab{cor-mg}
\end{Cor}
Moreover, we can
freely switch from $W_d$ in \rf{B.7} and $\widetilde W_d$ in \rf{Wtilde}.
Thus, by \rf{G.5.a} and \rf{G.5.b}, the higher-representation knot invariants
$W_d$ can be related to the fundamental-representation 
knot invariant $W_2$; besides, the limit
\beq
WW(C;\mu) := \lim_{d\to\infty,\lambda\to0, d\lambda = \mu}
\frac{W_d(C;\lambda)}d
\lab{WW}
\eeq
is well defined and, by \rf{D.23bistilde},
\beq
W_2(C;\lambda) = 2\frac d{d\lambda}[\lambda\ WW(C;\lambda)].
\lab{D.23bis}
\eeq
\vskip 10pt

We conclude this section with an important remark. Pure $BF$ theory is
known to be exact in saddle-point approximation \cite{BT2} as far as
the partition function is concerned. This result is a simple consequence
of the fact that one can arbitrarily change the ``Planck constant'' in front
of the action by simply rescaling the field $B$.  Since the partition
function (or the v.e.v.\ of an observable not containing $B$) is not
affected by this rescaling, one can send the Planck constant to
zero.

This, of course, cannot
be done when one computes the v.e.v.\ of an observable containing $B$, as
in \rf{B.7} or \rf{Wtilde}. In this case one sees that the parameter
$\lambda$ actually plays the role of the Planck constant, and of course,
$W_d$ is not independent of $\lambda$.  

Notice, however, 
that $WW$ corresponds to the limit $\lambda\to0$,
so it should be possible to compute it by using the saddle-point 
approximation. Then, by \rf{D.23bis}, it is possible to recover 
$W_2$---and hence, by \rf{G.5.a} and \rf{G.5.b}, all the $W_d$'s---from $WW$.

Thus, even if pure $BF$ theory with $B$-dependent observables is {\em not}
exact in saddle-point approximation, the saddle-point approximation 
turns out to be all that which one needs.

\section{Turning Off the Cosmological Constant}
\lab{sec-TOCC}
In this section we want to show that, in the limit of vanishing
cosmological,
the $BF$-theory generating function \rf{B.3} is related to both the
pure-$BF$-theory generating function \rf{B.7} and the Melvin--Morton
function \rf{MMf}.

To establish the former relation, we first observe that, in order
for the significant observable
$\gamma_1$ to survive in the limit $\kappa\to0$, we have to 
send $|x|\to\infty$ at the same time
with the prescription that $x\kappa=\lambda$ be finite.
If we work with the standard framing,
by Thm.\ \ref{th-sf} the observable $\gamma_0$ can be replaced
by the identity of the group; so \rf{B.3} diverges as $e^x$;
thus, we are led to consider
\beq
\Etilde(C;x,h):=e^{-x}\, E(C;x,h).
\lab{defEtilde}
\eeq
By using the exponential representation of \rf{B.3}, we can write
\[
\Etilde(C;x,h) = \aveBF{\trqw\exp\{x[\kappa\gamma_1(C,\xnot) 
+ O(\kappa^2)]\}
}\kappa
\]
and get
\beq
\lim_{|x|\to\infty,\ h\to0,\ xh=\hbar} \Etilde(C;x,h)
= W_2(C;\lambda),
\quad \hbar=4\pi i\lambda.
\lab{B.6}
\eeq
In App.\ \ref{app-kappa0}, we give a more careful proof of \rf{B.6}.

Notice that the limit in \rf{B.6} holds irrespectively of how $x$ 
is sent to infinity in the complex plane.
Therefore, if one sees $\Etilde$ as a meromorphic function of $x$ and $\hbar$, 
this implies that the regular part in $x$ of $\Etilde$ vanishes.  Thus,
$\Etilde$, now rewritten in terms of $x$ and $h$, does not contain more powers 
of $x$ than of $h$.  In other words, if $\Etilde$ is expanded as a rational
power series,
\beq
\Etilde(C;x,h) = \sum_{l,m=0}^\infty \epsilon_{lm}(C)\ x^lh^m,
\lab{B.4}
\eeq
then $\epsilon_\dodo$ is an upper triangular matrix,
\beq
\epsilon_{jm}(C) = 0, \quad\mbox{if $j>m$}.
\lab{epsilon}
\eeq

The relation between the $BF$-theory generating function and the
Melvin--Morton function is obtained by exploiting \rf{B.16} (s.\ App.\ 
\ref{app-c} for details).  We arrive to the following conclusions:
\beq
\epsilon_{\dodo}\ \mbox{upper triangular}
\Longleftrightarrow
b_{\dodo}\ \mbox{upper triangular},
\lab{epsb}
\eeq
and 
\beq
\epsilon_{nn}(C) = 2(1+n)\, b_{nn}(C),
\lab{cor2}
\eeq
where $\epsilon_{\dodo}$ and $b_{\dodo}$ are defined, respectively, 
in \rf{B.4} and \rf{expand}.
Notice that \rf{epsb} and \rf{epsilon} give a new proof of \rf{Con1MM}.

Since $E$ is related to $W_2$ and to $JJ$ in the same limit $\kappa\to0$,
we deduce that a relation exists between $W_2$ and $JJ$.  Actually,
by \rf{expand}, \rf{B.6}, \rf{B.4} and \rf{cor2}, we get
\beq
W_2(C;\lambda)=2\frac d{d\hbar}\left[{
\hbar\, JJ(C;\hbar)
}\right],
\quad \hbar = 4\pi i\lambda.
\lab{D.23}
\eeq
In order to recognize the knot invariant given by $W_2$, we can now 
resort to \rf{Con2MM} and obtain
\beq
W_2(C;\lambda)=2\left[{
1+\left({\frac z2}\right)^2
}\right]^{\frac 12}\ 
\frac d{dz} \frac z{\Delta(C;z)},
\quad z=2i\sin(2\pi\lambda),
\lab{D.26}
\eeq
and in particular,
\beq
W_2(\bigcirc;\lambda)=2\left[{
1+\left({\frac z2}\right)^2
}\right]^{\frac 12},
\quad z=2i\sin(2\pi\lambda),
\lab{D.27}
\eeq
where $\bigcirc$ is the unknot and we have chosen the normalization
$\Delta(\bigcirc;z)=1$, cfr.\ Sec.\ \ref{sec-CJF}.
Thus, the normalized knot invariant 
\beq
\ave C_2(\lambda) := \frac {W_2(C;\lambda)}{W_2(\bigcirc;\lambda)},
\lab{c2}
\eeq
defined in \cite{CCM}, satisfies 
\beq
\ave C_2(\lambda)= \frac d{dz} \frac z{\Delta(C;z)},
\quad z=2i\sin(2\pi\lambda).
\lab{c2bis}
\eeq

In conclusion, we have shown that $W_2$ is related to the first derivative
of the inverse of the Alexander--Conway polynomial.  
By \rf{G.5.a} and \rf{G.5.b}, $W_d$ is given by a finite sum of $W_2$'s,
evaluated at different $\lambda$'s.  Concerning the limit $d\to\infty$,
we see that \rf{D.23bis} and \rf{D.23}, together with the property
$WW(C;0)=JJ(C;0)=1$, imply
\beq
WW(C;\lambda) = JJ(C;\hbar),\quad \hbar=4\pi i\lambda,
\lab{WWJJ}
\eeq
or, because of \rf{Con2MM},
\beq
4\pi i\lambda\ WW(C;\lambda) = \frac z{\Delta(C;z)},
\quad z=2i\sin(2\pi\lambda).
\lab{WWDelta}
\eeq
The above results, together with Cor.\ \ref{cor-mg}, prove the following
\begin{Th}
The set of the unframed knot invariants that can be obtained from 
the $SU(2)$--$BF$ theory coincides with the set of the coefficents of
(the inverse of) the Alexander--Conway polynomial.
\lab{th-BF}
\end{Th}

In Ref.\ \cite{CCM}, the authors conjectured
a relation between the Alexander--Conway polynomyals and the pure $BF$ 
theory based on second-order calculations in the perturbative expansion.
Thm.\ \ref{th-BF} supersedes this conjecture and provides
the correct relation.

\section{The Perturbative Expansion of $W_d$}
\lab{sec-PEW}
In this section we want to compare the results we have proved in 
Sec.\ \ref{sec-TOCC} with those obtained in \cite{CCM}
in the framework of ``perturbative'' $BF$ theory.

By ``perturbative'' evaluation of \rf{B.7},
one means an expansion of
$W_d(C;\lambda)$ (or of its generalization $W_{G,R}(C;\lambda)$,
where $G$ is a compact group and $R$ a 
representation) in powers of $\lambda$,
\beq
W_{G,R}(C;\lambda) = \sum_{n=0}^\infty w_n(G,R;C)\ \lambda^n,
\lab{Wexpand}
\eeq
where the knot invariants $w_n(G,R;C)$ (actually, they are Vassiliev invariants
\cite{Vass}) are computed in terms of Feynman integrals.

The first property of the expansion \rf{Wexpand} shown in \cite{CCM} is
that, owing to a symmetry of the corresponding Feynman integrals, odd order 
terms vanish,
\beq
w_{2n+1}(G,R;C) = 0.
\lab{wodd}
\eeq
Thus, $W$ is an even function
\beq
W_{G,R}(C;-\lambda) = W_{G,R}(C;\lambda).
\lab{Weven}
\eeq
In the case $G=SU(2)$, by using \rf{D.26}, \rf{G.5.a} and \rf{G.5.b},
we see that \rf{Weven} is
in accordance with the fact that the Alexander--Conway
polynomial $\Delta(C;z)$---for a single knot $C$---is an even function of $z$. 

A further computation done in \cite{CCM}, for the case $G=SU(N)$,
showed that, up to the second order, \rf{Wexpand} reads
\beq
W_{G,R}(C;\lambda) = \dim R\ \left[{
1 + (4\pi\lambda)^2\ c_2(R)\, c_v\ \rho(C) + O(\lambda^4)
}\right],
\lab{E.6}
\eeq
where:
\begin{itemize}
\item $\dim R$ and $c_2(R)$ are, respectively, the dimension and the 
quadratic Casimir of the representation $R$
\item $c_v=N$ is the quadratic Casimir of the adjoint representation
\item  $\rho(C)$ is the knot invariant studied in \cite{GMM,BN-th}
in the framework of the Chern--Simons theory
\end{itemize}
In \cite{GMM} and \cite{BN-th}, $\rho(C)$ was proved to be related
to the second coefficient of the Alexander--Conway polynomial; viz.,
if one writes
\beq
\Delta(C;z) = \sum_{n=0}^\infty a_n(C)\, z^n,
\quad a_0(C)=1,\ a_1(C)=0
\lab{Deltaexpand}
\eeq
(where only a finite number of coefficients are nonvanishing), then
\beq
\rho(C) = 2\, a_2(C) + \rho(\bigcirc).
\lab{rhoC}
\eeq
Moreover, a direct computation \cite{GMM} shows that
\beq
\rho(\bigcirc) = -\frac 1{12}.
\lab{rhoO}
\eeq
Thus, one can write \rf{E.6} as
\beq
W_{G,R}(C;\lambda) = \dim R\ \left[{
1 + (4\pi\lambda)^2\ 2\, c_2(R)\, c_v \left({
a_2(C)-\frac1{24}
}\right) + O(\lambda^4)
}\right].
\lab{E.6'}
\eeq
If we now expand \rf{D.26} as
\begin{eqnarray*}
W_2(C;\lambda) &=& 2 \left[{
1 + \frac18z^2 + O(z^4)
}\right] \left[{
1 - 3a_2(C) z^2 + O(z^4)
}\right] =\\
 &=& 2 \left[{
1 + 3z^2  \left({
\frac1{24}-a_2(C)
}\right)+O(z^4)
}\right],\quad z = 4\pi i\lambda + O(\lambda^2),
\end{eqnarray*}
we get a complete
agreement with \rf{E.6'} in the case where $G=SU(2)$ and the fundamental
representation, $R_2$, is chosen; for, in this case, $c_v=2$ and
$c_2(R_2) = 3/4$.

We can also compare the second-order expansion of \rf{G.5.a} and
\rf{G.5.b} with \rf{E.6'} in the case where $G=SU(2)$ and $R_d$ is
the irreducible representation of dimension $d$.
This is easily done by noticing that one can easily compute the sums
appearing in \rf{G.5.a} and \rf{G.5.b} as
\begin{eqnarray*}
\sum_{l=1}^k (2l-1)^2 &=& \frac13 k(4k^2-1),\\
\sum_{l=1}^k l^2 &=& \frac23 k(k+1)(2k+1). 
\end{eqnarray*}
Thus, one achieves complete agreement with \rf{E.6'} since
\beqy
c_2(R_{2k}) &=& \frac{(2k-1)(2k+1)}4,\lab{c2even}\\
c_2(R_{2k+1}) &=& k(k+1). \lab{c2odd}
\eeqy

By using \rf{G.5.a} and \rf{G.5.b}, it is possible to
see that---in agreement with the perturbative result of \cite{CCM}---the 
coefficients $w_n$ in \rf{Wexpand} are given by the product
of a function depending only on the representation and a function
depending only on the knot, viz.,
\beq
w_{2n}(SU(2),R_d;C) = d\, g_{2n}(d)\, v_{2n}(C).
\lab{wgv}
\eeq
If we normalize $g_{2n}(2)=1/2$, then, by \rf{G.5.a} and \rf{G.5.b},
we get an explicit formula for the $g$'s,
\beqy
g_{2n}(2k) &=& \frac1{2k}\sum_{l=1}^k (2l-1)^{2n},\\
g_{2n}(2k+1) &=& \frac1{2k+1}\sum_{l=1}^k (2l)^{2n}.
\eeqy
By \rf{c2even} and \rf{c2odd}, we can also express the $g$'s
in terms of the quadratic Casimirs:
\beq
\begin{array}{rcl}
g_0 &=& \frac12,\\
g_2 &=& \frac13 c_2 c_v,\\
g_4 &=& \frac1{15}\left[{
6(c_2c_v)^2 - c_2c_v^3
}\right],\\
g_6 &=& \frac1{21}\left[{
12(c_2c_v)^3 - 6 c_2^2c_v^4 + c_2c_v^5
}\right],\\
\cdots
\end{array}
\lab{g's}
\eeq
In general, the highest power of $c_2$ in $g_{2n}$ is $n$.

In \cite{CCM}, it was conjectured that the {\em only} group factor
to appear in $g_{2n}$ were $(c_2c_v)^n$. 
{}From this conjecture, together with some field-theoretical arguments, 
it was concluded that, for some $t$,
\beq
\frac{W_{G,R}(C,\lambda)}{W_{G,R}(\bigcirc,\lambda)} = [\Delta(C;z)]^t,
\quad z\propto \sqrt{c_2(R)\, c_v}\,\lambda + O(\lambda^3).
\lab{conjCCM}
\eeq
By \rf{g's}, we see that the conjecture is wrong;
by \rf{D.23}, we see that so is the conclusion \rf{conjCCM}.  

It is however interesting
to show that, if we retain only the terms $(c_2c_v)^n$ in 
\rf{g's}, then the corresponding knot 
invariant
is actually the $(-1)$-power of the Alexander--Conway 
polynomial. More precisely, we define the ``truncated'' coefficients
\beq
g_{2n}^{\rm (tr)}(d) =  [c_2(R_d)\, c_v]^n\, \gamma_{2n}
\lab{gtr}
\eeq
as the coefficients obtained by neglecting lower powers
of $c_2(R_d)$ in \rf{g's}, and the ``truncated'' knot invariant as
\beq
W_d^{\rm (tr)}(C;\lambda) := \sum_{n=0}^\infty d\, g_{2n}^{\rm (tr)}(d)
\,v_{2n}(C)\, \lambda^{2n}.
\lab{Wtr}
\eeq
Since $c_v=2$ and $c_2(R_d)\sim d^2/4$ as $d\to\infty$, the truncated
coefficients $g_{2n}^{\rm (tr)}(d)$ are the leading terms of the
true coefficients $g_{2n}(d)$; thus, by \rf{WW} and \rf{gtr},
\beq
WW(C;\lambda) = \sum_{n=0}^\infty \frac{\gamma_{2n}}{2^n}\, v_{2n}(C) 
\,\lambda^{2n}.
\lab{WWtr}
\eeq
By comparing \rf{WWtr} with \rf{Wtr}, we eventually get
\beq
W_d^{\rm (tr)}(C;\lambda)=d\ WW(C;\sqrt{2\,c_2(R_d)\, c_v}\,\lambda);
\eeq
thus, by \rf{WWDelta}, we see that
\beq
\frac{W_d^{\rm (tr)}(C;\lambda)}{W_d^{\rm (tr)}(\bigcirc;\lambda)}
= \frac 1{\Delta(C;z)},\quad z = 2i\sin\left({
2\pi \sqrt{2\,c_2(R_d)\, c_v}\,\lambda
}\right).
\lab{WtrDelta}
\eeq
Thus, if the truncated knot invariants are used, \rf{conjCCM} holds.

\section{Conclusions}
In this paper we have discussed the unframed knot invariants coming 
from $BF$ theories. Even if most of our results hold only for $SU(2)$,
we point out that Thm.\ \ref{th-sf} 
and Cor.\ \ref{cor-mg}, as well as the computation in 
App.\ \ref{app-kappa0} (with a slight abuse of notation), hold in general.

It would be interesting to generalize some of the other results to different
groups and to consider links as well.
However, the present case---i.e., knot observables in the theory based
on $SU(2)$---seems to be interesting enough
to deserve further investigation.

Indeed, as we have noticed at the
end of Sec.\ \ref{sec-PBF}, it turns out that a saddle-point computation
is enough to completely describe the pure $BF$ theory. We defer
to a forthcoming paper the related functional-integral computation. 

Notice that this property of pure $BF$ theory sets it
at the boundary between TQFTs of Witten's and Schwarz's type: The former
are twisted supersymmetric gauge theories (s., e.g., \cite{WittSUSY}),
whose main property---an effect of the twisted supersymmetry---is their
independence of both the metric and the coupling constant, 
which makes them topological as well as exact in saddle-point approximation.
The latter are topological gauge theories, as the Chern--Simons or the
$BF$ theories, whose dependence on the coupling constant is unavoidable.
The pure $BF$ theory formally belongs to the latter type, but, as the
theories of the former type, is completely determined by its
weak-coupling limit.

Moreover, as we have proved in Sec.\ \ref{sec-TOCC}, 
the pure $BF$ theory corresponds
to the first diagonal in the $(h,d)$ expansion of the colored Jones function.
A description of the upper diagonals is still missing (s.\ \cite{Roz2}
for a first attempt); so it is natural to look for generalizations
of the pure $BF$ theory, i.e., for further variations of the Chern--Simons
theory, that could correspond to these upper diagonals,
and possibly, give them a
better understanding (s.\ \cite{KSA} for a different approach). 
We are investigating along these lines.

\section*{Acknowledgements}
This work was supported by INFN Grant No.\ 5077/94.

I thank A.~Jaffe, T.~Kerler, A.~Le\'sniewski, M.~Martellini and M.~Rinaldi
for helpful conversations. I am especially thankful to P.~Cotta-Ramusino
for a number of very useful discussions and for constant advice.

\appendix
\section{Some Useful Identities for the Characters of $SU(2)$}
\lab{app-A}
In this section we recall some properties of the group
$SU(2)$ that are necessary for this paper.  In particular,
we are interested in an identity relating
$\trqw g^n$
to the traces of $g$ in other representations, or in other words,
we want to give $\trqw g^n$ an expansion in primitive characters,
\beq
\trqw g^n = \sum_{k=0}^\infty \phi_{nk}\, \chi_{\frac k2}(g),
\lab{A.2}
\eeq
where
\beq
\chi_s(g) := \tr s g, \quad 2s\in{\bf Z}.
\lab{chi}
\eeq
Owing to the Peter--Weil Theorem, this can be done
since $\trqw g^n$ depends only on the conjugacy class of $g$ and belongs to 
$L^2(SU(2),\complessi)$.  By a conjugacy transformation, one can always write
$g$ as
\beq
g = h\, e^{i\alpha R_3}\, h^\dagger,
\lab{A.3}
\eeq
where $R_3$ is the (hermitean) generator of the Cartan subalgebra.
Noticing that
the spectrum of $R_3$ in the representation of spin $l$ is given by
$\{-2l, -2(l-1),\ldots, 2(l-1), 2l\}$, one gets
\beq
\trqw g^n = e^{in\alpha}+e^{-in\alpha},
\lab{A.4}
\eeq
and
\beq
\tr l g = e^{2il\alpha} + e^{2i(l-1)\alpha} + \cdots
+ e^{-2il\alpha} = \sum_{m=-l}^l e^{2im\alpha},
\lab{A.5}
\eeq
where the last sum is meant to be over integers or half-integers according
to the fact that $l$ is integer or half-integer.  
By \rf{A.4} and \rf{A.5}, it immediately follows that, for $n\ge2$,
\beq
\trqw g^n = \tr{\frac n2}g - \tr{\frac n2 - 1}g, \quad n\in{\bf Z}.
\lab{A.7}
\eeq
By induction one can also prove that
\beq
\tr s g = \left\{
\begin{array}{ll}
\sum_{l=1}^{s+1/2} \trqw g^{2l-1} & \mbox{if $s$ is half-integer}\\
1+\sum_{l=1}^s \trqw g^{2l} & \mbox{if $s$ is integer}
\end{array}
\right.
\lab{gn}
\eeq
Moreover, by \rf{gn} and by the fact that
\[
\trqw \exp(i {\bf a\cdot R}) = 2\cos a,
\]
with ${\bf a}\in\reali^3$ and ${\bf R}=(R^1,R^2,R^3)$ the Pauli matrices,
one also gets
\beq
\tr s \exp(i {\bf a\cdot R}) = \frac{\sin(da)}{\sin(a)},
\lab{trs}
\eeq
where $d=2s+1$ is the dimension of the representation of spin $s$.

Since, for a given $g$, $\tr {\frac n2} g$ is a map defined on the positive
integers, its analytic continuation over the whole complex plane is
uniquely defined and is actually given by \rf{trs}.

\section{Properties of the Generating Function for Colored Jones Functions}
\lab{app-f}
In this section we exploit the properties of 
the generating function for colored Jones functions defined in \rf{gf}.
By \rf{expand}, we can write
\beq
f(C;x,h) = \sum_{d=1}^\infty \frac{x^d}{d!}\, d\ 
\sum_{j,m=0}^{\infty} b_{jm}(C)\ (d-1)^jh^m=
x\ \sum_{j,m=0}^{\infty} B_j(x)\ b_{jm}(C)\ h^m,
\lab{fexpand}
\eeq
where
\beq
B_j(x) = \sum_{d=0}^\infty \frac{x^d}{d!} d^j =
\left({\frac{\de^j}{\de\alpha^j}}\right)_{\big|_{\alpha=\log x}}
\sum_{d=0}^\infty \frac{e^{\alpha d}}{d!} =
\left({\frac{\de^j}{\de\alpha^j}}\right)_{\big|_{\alpha=\log x}} 
\!\!\!\!\!\!\!\!\!\!\!\!    \exp(e^\alpha).
\lab{App.2}
\eeq
By repeatedly applying Leibniz's rule, we get
\beq
B_j(x) = e^x P_j(x) = e^x [x^j + O(x^{j-1})],
\lab{bj}
\eeq
where $P_j(x)$ is a polynomial of degree $j$ starting with $x^j$, viz.,
\beq
P_j(x) = \sum_{l=0}^j x^l\, c_{lj}, \quad c_{jj}=1.
\lab{Pj}
\eeq
Notice that \rf{Pj} implies that $c_\dodo$ is an upper triangular matrix:
\beq
c_{lj}(C) = 0, \quad\mbox{if $l>j$}.
\lab{c}
\eeq

The coefficients $c_{lj}$ can be easily computed if we consider the
following generating function
\beq
B(\rho,x) := \sum_{j=0}^\infty \frac {\rho^j}{j!}\, B_j(x),
\lab{App.7}
\eeq
which, by \rf{App.2}, can be written as
\beq
B(\rho,x) = \exp{e^{(\alpha+\rho)}}_{\big|_{\alpha=\log x}} =
e^x\, \exp[x\,(e^\rho-1)].
\lab{App.8}
\eeq
Thus, by \rf{bj}, we obtain
\beq
P(\rho,x) := \sum_{j=0}^\infty \frac {\rho^j}{j!}\, P_j(x)=
\sum_{n=0}^\infty \frac{x^n}{n!}\ (e^\rho-1)^n,
\lab{App.9}
\eeq
while, by \rf{Pj}, we have
\beq
P(\rho,x) = \sum_{n,j=0}^\infty x^n\, c_{nj}\,  \frac {\rho^j}{j!}.
\lab{App.9'}
\eeq
By comparing the two different expansions of $P$ in powers of $\rho$
given by \rf{App.9} and \rf{App.9'}, we eventually get
\beq
c_{nj} = \frac {(-1)^n}{n!} \sum_{l=0}^n \bino nl\, (-1)^l\, l^j,
\lab{App.10}
\eeq
where, by convention, $0^0=1$.

We want now to compute the coefficients $\btilde_\dodo$ defined in
\rf{C.18}. By \rf{C.13}, \rf{fexpand} and \rf{bj}, we can write
\beq
\ftilde(C;x,h) = \sum_{j,m=0}^\infty P_j(x)\ b_{jm}(C)\ h^m;
\lab{C.18'}
\eeq
so by \rf{C.18}, 
\beq
\btilde_{nm}(C) = \sum_{j=0}^\infty c_{nj}\ b_{jm}(C).
\lab{C.19}
\eeq

Since the matrices $\btilde_\dodo$ and $b_\dodo$ are related by the
matrix $c_\dodo$ that, by \rf{c}, is upper triangular, we conclude that
\rf{bbtilde} holds.  Moreover, if one knows that either $\btilde_\dodo$
or $b_\dodo$ is upper triangular---and by \rf{Con1MM} we know this to
obtain for the latter---then one has
\beq
\btilde_{mm}(C) = b_{mm}(C).
\lab{tildebb}
\eeq
This implies that, in the limit considered in \rf{C.14}, we can write
\[
\ftilde(C;x,h) = \sum_{m=0}^\infty \btilde_{mm}(C)\ \hbar^m + 
O\left(\frac\hbar x\right)=
\sum_{m=0}^\infty b_{mm}(C)\ \hbar^m + 
O\left(\frac\hbar x\right);
\]
so by \rf{MMf}, \rf{C.14} holds.  In order to prove \rf{C.16}, we have
only to notice that the operator $(x\, \de/\de x)^n$, when applied to $x^j$,
simply produces a factor $j^n$; so
\[
\left({x\frac\de{\de x}}\right)^n
\ftilde(C;x,h) = \sum_{m=0}^\infty m^n\,\btilde_{mm}(C)\ \hbar^m + 
O\left(\frac\hbar x\right)=
\sum_{m=0}^\infty m^n\,b_{mm}(C)\ \hbar^m + 
O\left(\frac\hbar x\right).
\]
By comparison with the effect of the operator $(\hbar\,d/d\hbar)^n$ on
\rf{MMf}, we get \rf{C.16}. 

\section{From $E$ to $W_2$ as $\kappa\to0$}
\lab{app-kappa0}
In this section we want to clarify the limit in \rf{B.6}. By
\rf{defEtilde}, \rf{B.3} and \rf{Gamma}, we have
\beq
\Etilde(C;x,h)=e^{-x}\sum_{n=0}^\infty \frac{x^n}{n!}
\aveBF{\trqw[\gamma_0(C,\xnot) + \kappa\gamma_1(C,\xnot) + O(\kappa^2)]^n
}\kappa,
\lab{k1}
\eeq
where $h=4\pi i\kappa$.

In Reff.\ \cite{C,CCFM}, it is shown that, in an observable, the field $B$ 
represents a source for $F+\kappa^2B\wedge B$, while the field $A$ is a source
for $d_AB$ only.  Since a variation in the framing (i.e., the companion
knot along which we integrate the field $A$) is still given by inserting
a curvature term,
we can repeat the steps of the proof 
of Thm.\ \ref{th-sf} and show that, if the standard framing is chosen,
$\Gamma_0$ can be replaced by $I+O(\kappa^2)$, where $I$ is the group 
identity. As a consequence, \rf{k1} now reads
\beq
\Etilde(C;x,h)=e^{-x}\sum_{n=0}^\infty \frac{x^n}{n!}
\aveBF{\trqw[I + \kappa\gamma_1(C,\xnot) + O(\kappa^2)]^n
}\kappa.
\lab{k2}
\eeq
By using Newton's binomial formula, one can easily prove that \rf{k2} can
also be written as
\[
\Etilde(C;x,h)=\sum_{n=0}^\infty \frac{x^n}{n!}
\aveBF{\trqw[\kappa\gamma_1(C,\xnot) + O(\kappa^2)]^n
}\kappa,
\]
or setting $\lambda=x\kappa$,
\beq
\Etilde(C;x,h)=\sum_{n=0}^\infty \frac{\lambda^n}{n!}
\aveBF{\trqw[\gamma_1(C,\xnot) + O(\kappa)]^n
}\kappa.
\lab{k3}
\eeq
Now, sending $\kappa\to0$, with $\lambda$ fixed, gives \rf{B.6}.

\section{From $E$ to $JJ$ as $h\to0$}
\lab{app-c}
In this section we consider the relation between the $BF$-theory
generating function $E$ and the Melvin--Morton function as the expansion
parameter $h$ is sent to zero. 

The starting point is relation \rf{B.16}, which, by \rf{C.13}, can be
rewritten in terms of $\tilde f$.
Actually, the second
term on the r.h.s. of \rf{B.16} is easily seen to be
\beq
\de_xf(C;x,h) = e^x\,(x+1)\,\ftilde(C;x,h)+
e^x\, x\,\de_x\ftilde(C;x,h),
\lab{D.1}
\eeq
while the computation of the last term requires a piece more of work. 
First of all, it is useful to introduce the following notation
\beq
I_n[g](x) = \int_0^x d\xi\ e^\xi\, g^{(n)}(\xi),
\lab{D.2}
\eeq  
where $g$ is a generic analytic function, and $g^{(n)}$ its $n$th derivative. 
By integrating by parts, one can 
prove the recursion rule
\beq
I_n[g](x) = e^x\,  g^{(n)}(x) - g^{(n)}(0) - I_{n+1}[g](x),
\lab{D.3}
\eeq
which implies that
\beq
I_0[g](x) = e^x\sum_{n=0}^\infty (-1)^n\, g^{(n)}(x)
- \sum_{n=0}^\infty (-1)^n\, g^{(n)}(0).
\lab{D.4}
\eeq
Now we set
\beq
g(x) = x\, \ftilde(C;x,h)
\lab{D.5}
\eeq
(where $C$ and $h$ are supposed to be fixed), and note that, accordingly,
\beq
g^{(n)}(x) = n\,\de_x^{n-1} \ftilde(C;x,h) + x\, \de_x^{n} \ftilde(C;x,h).
\lab{D.16}
\eeq
Then we use \rf{D.4} to compute
the last term in \rf{B.16} as
\beq
\int_0^x d\xi\ f(C;\xi,h) =
e^x\left[{
(x-1)\,\ftilde(C;x,h) - x\,\de_x\ftilde(C;x,h) + R(C;x,h)
}\right]
+c(C;h),
\lab{D.19}
\eeq
where
\beq
R(C;x,h) = \sum_{n=2}^\infty (-1)^n
\left[{
n\,\de_x^{n-1} \ftilde(C;x,h) + x\, \de_x^{n} \ftilde(C;x,h)
}\right]
\lab{D.17}
\eeq
and
\beq
c(C;h) = \sum_{n=0}^\infty (-1)^n\, (n+1)\, \de_x^n \ftilde(C;0,h).
\lab{D.20}
\eeq
Therefore, by \rf{defEtilde}, \rf{B.16}, \rf{D.1} and \rf{D.19}, we get
\beq
\Etilde(C;x,h) = 2
\left[{
\ftilde(C;x,h) + x\, \de_x \ftilde(C;x,h)
}\right]
+R(C;x,h) + 
e^{-x}\left[{1 - c(C;h)}\right].
\lab{D.21}
\eeq
We defer to the end of this section the proof that $c(C;h)=1$, cfr.\ \rf{F.17};
as a consequence of this fact, \rf{D.21} actually reads 
\beq
\Etilde(C;x,h) = 2
\left[{
\ftilde(C;x,h) + x\, \de_x \ftilde(C;x,h)
}\right]
+R(C;x,h).
\lab{D.21'}
\eeq
Now we want to re-express \rf{D.21'} as a relation between the coefficients
of the power series expansions of $\Etilde$ and $\ftilde$.  Indeed, by
\rf{C.18}, the terms in square brackets in \rf{D.21'} can be written as
\[
\sum_{n,l=0}^\infty (1+l)\ \btilde_{lm}(C)\ x^l h^m,
\]
while
\beq
R(C;x,h) = \sum_{n=2}^\infty\ \sum_{l,m=0}^\infty
(-1)^n\, \frac{(n+l)!}{l!}\ \btilde_{l+n-1,m}(C)\ x^l h^m.
\lab{F.3}
\eeq
Therefore, by \rf{B.4}, we eventually get
\beq
\epsilon_{lm}(C) = \sum_{j=0}^\infty \eta_{lj}\, \btilde_{jm}(C),
\lab{F.4}
\eeq
with
\beq
\eta_{lj} = 2\, (1+l)\, \delta_{lj} + \sum_{n=2}^\infty
(-1)^n\, \frac{(n+l)!}{l!}\ \delta_{l+n-1,j}
\lab{F.5}
\eeq
($\delta_\dodo$ being the Kronecker delta).  
Since $\eta_\dodo$ is an upper
triangular matrix, by \rf{F.4} we have that
\beq
\epsilon_{\dodo}\ \mbox{upper triangular}
\Longleftrightarrow
\btilde_{\dodo}\ \mbox{upper triangular}.
\lab{epsbtilde}
\eeq
Finally, by \rf{bbtilde}, we conclude that \rf{epsb} holds.
Moreover, \rf{F.5} and \rf{tildebb} imply \rf{cor2}.
\vskip 10pt

We conclude this section by showing that the function $c(C;h)$, 
defined in \rf{D.20}, is a constant equal to one.  

We start by considering the expansion of $c$ in powers of $h$. By \rf{C.18},
we get
\beq
c(C;h) = \sum_{n,m=0}^\infty (-1)^n\, (n+1)!\ \btilde_{nm}(C)\ h^m =
\sum_{j,m=0}^\infty D_j\ b_{jm}(C)\ h^m,
\lab{F.8}
\eeq
where, by \rf{C.19}, 
\beq
D_j = \sum_{n=0}^\infty (-1)^n\, (n+1)!\ c_{nj}.
\lab{F.13}
\eeq
Notice that, by \rf{c}, \rf{F.13} is actually a finite sum.
By \rf{App.10}, we can also write
\beq
D_j = \sum_{n=0}^\infty\,\sum_{l=0}^n\,
(n+1)\, \bino nl\, (-1)^l\, l^j.
\lab{F.14}
\eeq
Let us consider now the generating function
\beq
D(\rho):=\sum_{j=0}^\infty \frac {\rho^j}{j!} D_j.
\lab{F.D}
\eeq
By \rf{F.14}, it follows that
\beq
D(\rho) = \sum_{n=0}^\infty (n+1)\ 
\sum_{l=0}^n \bino nl\, (-1)^l\, e^{\rho l} =
\sum_{n=0}^\infty (n+1)\, (1-e^\rho)^n = e^{-2\rho}.
\lab{F.21}
\eeq
By comparing the expansion of $D$ in powers of $\rho$ in \rf{F.D}
with the expansion of $e^{-2\rho}$, we conclude that 
\beq
D_j=(-2)^j.
\lab{F.14.2}
\eeq
Thus, \rf{F.8} reads
\beq
c(C;h)=\sum_{j=0}^\infty b_{jm}(C)\ (-2)^j\,h^m =
-J_{-1}(C;h),
\lab{F.16}
\eeq
where the last identity follows from \rf{expand}.
However, by \rf{Jodd}, we have
\[
J_{-1}(C;h) = -J_1(C;h) = -1.
\] 
Therefore, we conclude that
\beq
c(C;h)=1.
\lab{F.17}
\eeq

\end{document}